\documentclass[12pt,a4paper,reqno]{article}

\usepackage[latin1]{inputenc} 
\usepackage[T1]{fontenc}
\usepackage{amscd}
\usepackage{amsmath}
\usepackage{amsfonts}
\usepackage{amssymb}
\usepackage{amsthm}
\usepackage{comment}
\usepackage[frenchb,british]{babel} 
\usepackage[a4paper,portrait]{geometry}
\ifx\pdfoutput\undefined
\usepackage[dvips,a4paper,plainpages=false]{hyperref}
\else
\pdfoutput=1\pdfcompresslevel=9
\usepackage[pdftex,a4paper,plainpages=false]{hyperref}
\fi
\hypersetup{colorlinks=true}
\hypersetup{pagecolor=black,linkcolor=black,citecolor=black,urlcolor=black}

\newcommand{\THMEN}{%
\theoremstyle{plain}%
\newtheorem{thm}{Theorem}[section]%
\newtheorem{cor}[thm]{Corollary}%
\newtheorem{prop}[thm]{Proposition}%
\newtheorem{lem}[thm]{Lemma}%
\theoremstyle{definition}%
\newtheorem{defi}[thm]{Definition}%
\theoremstyle{remark}%
\newtheorem{rem}[thm]{Remark}%
\newtheorem{xpl}[thm]{Example}%
\newtheorem{hyp}[thm]{Hypothesis}%
}

\newcommand{\EM}{\ensuremath}

\newcommand{\dE}{\EM{\mathbb{E}}}

\newcommand{\dL}{\EM{\mathbb{L}}}

\newcommand{\dN}{\EM{\mathbb{N}}}

\newcommand{\dP}{\EM{\mathbb{P}}}

\newcommand{\dR}{\EM{\mathbb{R}}}

\newcommand{\cC}{\EM{\mathcal{C}}}
\newcommand{\cD}{\EM{\mathcal{D}}}

\newcommand{\cF}{\EM{\mathcal{F}}}

\newcommand{\cH}{\EM{\mathcal{H}}}

\newcommand{\cM}{\EM{\mathcal{M}}}
\newcommand{\cN}{\EM{\mathcal{N}}}
\newcommand{\cO}{\EM{\mathcal{O}}}
\newcommand{\cP}{\EM{\mathcal{P}}}

\DeclareMathOperator{\Tr}{Tr}

\DeclareMathOperator{\ide}{Id}

\DeclareMathOperator{\ree}{Re}

\newcommand{\al}{\alpha}

\newcommand{\de}{\delta}

\newcommand{\Ga}{\Gamma}

\newcommand{\la}{\lambda}

\newcommand{\Om}{\Omega}

\newcommand{\Si}{\Sigma}

\newcommand{\te}{\theta}

\newcommand{\eps}{\varepsilon}
\newcommand{\ffi}{\varphi}

\newcommand{\lam}{\lambda_{\text{max}}^N}
\newcommand{\mum}{\mu_{\text{max}}^{M,N}}
\newcommand{\gue}{\text{GUE}(N)}
\newcommand{\lue}{\text{LUE}(N,M)}

\newcommand{\ot}{\otimes}
\newcommand{\od}{\odot}

\newcommand{\lsup}{\limsup}

\newcommand{\II}{1 \! \mathrm{I}}

\THMEN

\author{Yan Doumerc}
\date{University of Toulouse}
\title{An asymptotic link between LUE and GUE and its spectral interpretation}

\begin{document}
\maketitle 
\begin{abstract}
We use a matrix central-limit theorem which makes
the Gaussian Unitary Ensemble appear as a limit of the Laguerre
Unitary Ensemble together with an observation due to Johansson in order to
derive new representations for the eigenvalues of GUE. For instance, it is
possible to recover the celebrated equality in distribution between the
maximal eigenvalue of GUE and a last-passage time in some directed brownian
percolation. Similar identities for the other eigenvalues of GUE also appear.
\end{abstract}

\section{Introduction}

The most famous ensembles of Hermitian random matrices are undoubtedly
the Gaussian Unitary Ensemble (GUE) and the Laguerre Unitary Ensemble
(LUE). Let $(X_{i,j})_{1\le i<j\le N}$ (respectively $(X_{i,i})_{1\le
  i\le N}$) be complex (respectively real) standard independent Gaussian
variables ($\dE(X_{i,j})=0$, $\dE(|X_{i,j}|^2)=1$) and let $X_{i,j}=\overline{X_{j,i}}$ for $i>j$. The $\gue$ is
defined to be the random matrix $X^N=(X_{i,j})_{1\le i,j \le N}$. It induces
the following probability measure on the space $\cH_N$ of $N \times N$ Hermitian matrices: 
\begin{equation}
\label{densg}
P_N(dH)=Z_N^{-1} \exp\Big(-\frac{1}{2}
\Tr(H^2)\Big)dH
\end{equation}
where $dH$ is Lebesgue measure on $\cH_N$. In the same
way, if $M\ge N$ and $A^{N,M}$ is a $N \times M$ matrix whose entries are complex standard independent Gaussian variables, then $\lue$ is
defined to be the random $N \times N$ matrix $Y^{N,M}=A^{N,M} (A^{N,M})^{\ast}$ where
$\ast$ stands for the conjugate of the transposed matrix. Alternatively, $\lue$
corresponds to the following measure on $\cH_N$: 
\begin{equation}
\label{densl}
P_{N,M}(dH)=Z_{N,M}^{-1}
(\det H)^{M-N} \exp(-\Tr H) \II_{H \ge 0}dH \; .
\end{equation}
Here we give a proof of the fact that $\gue$ is the limit in distribution of
$\lue$ as $M \to \infty$ in the following asymptotic regime: 
\begin{equation}
\label{limite}
\frac{Y^{N,M}-M\ide_N}{\sqrt{M}}\overset{d}{\underset{M \to \infty}{\longrightarrow}} X^N \; .
\end{equation}
We also prove such a central-limit theorem at a process level when the
Gaussian entries of the matrices are replaced by Brownian motions. The
convergence takes place for the trajectories of the eigenvalues. 

Next, we make use of this matrix central-limit theorem together with an
observation due to Johansson \cite{johansson-shape} and an invariance principle for a
last-passsage time due to Glynn and Whitt \cite{gly-whi} in order to recover the
following celebrated equality in distribution between the maximal eigenvalue
$\lam$ of $\gue$ and some functional of standard $N$-dimensional Brownian
motion $(B_i)_{1\le i\le N}$ as
\begin{equation}
\label{qqq}
\lam \overset{d}{=} \sup_{0=t_0\le\cdots\le t_N =1} \sum\limits_{i=1}^{N}
(B_i(t_i)-B_i(t_{i-1})) \; .
\end{equation}
The right-hand side of (\ref{qqq}) can be thought of as a last-passage time in an oriented
Brownian percolation. Its discrete analogue for an oriented percolation
on the sites of $\dN^2$ is the object of Johansson's remark. The identity
(\ref{qqq}) first appeared in \cite{baryshnikov} and \cite{gra-tra-wid}. Very recently,
O'Connell and Yor shed a remarkable light on this result in \cite{oconnell-representation}. Their
work involves a representation similar to (\ref{qqq}) for all the eigenvalues of $\gue$. We notice here that analogous formula can be written for all the eigenvalues of
$\lue$. On the one hand, seeing the particular expression of these
formula, a central-limit theorem can be established for them and the
limit variable $\Om$ is identified in terms of Brownian functionals. On the
other hand, the previous formulas for eigenvalues of $\lue$ converge, in
the limit given by (\ref{limite}), to the representation found in
\cite{oconnell-representation} for $\gue$ in terms of some path-transformation $\Ga$ of
Brownian motion. It is not immediately obvious to us that functionals
$\Ga$ and $\Om$ coincide. In particular, is this identity true pathwise or only in distribution?  

The matrix central-limit theorem is presented in Section \ref{clt} and its proof is postponed to the last section. In section \ref{int}, we described the consequences to eigenvalues representations and the connection with the O'Connell-Yor approach.

\section{The central-limit theorem}
\label{clt}

We start with the basic form of the central-limit theorem we investigate here.

\begin{thm}
\label{lim}
Let $Y^{N,M}$ and $X^N$ be taken respectively from $\lue$ and
$\gue$. Then 
\begin{equation}
\label{simple}
\frac{Y^{N,M}-M\ide_N}{\sqrt{M}}\overset{d}{\underset{M
    \to \infty}{\longrightarrow}} X^N \, .
\end{equation}
\end{thm}

\begin{rem}
The fact that each entry of the left-hand side in (\ref{simple})
converges to the corresponding entry of the right-hand side is just
the classical central-limit theorem. But (\ref{simple}) asserts that
there is also an asymptotic independence of the upper-diagonal
entries, which is not obvious from the explicit formula of entries.
\end{rem}

We turn to the process version of the previous result. Let
$A^{N,M}=(A_{i,j})$ be a $N \times M$ matrix whose entries are independent standard complex
Brownian motions. The Laguerre process is defined to be
$Y^{N,M}=A^{N,M}(A^{N,M})^{\ast}$. It is built in exactly the same way as
$\lue$ but with Brownian motions instead of Gaussian
variables. Similarly, we can define the Hermitian Brownian motion
$X^N$ as the process extension of $\gue$.  

\begin{thm}
\label{processus}
If $Y^{N,M}$ is the Laguerre process and $(X^N(t))_{t\ge 0}$ is Hermitian Brownian motion, then:
\begin{equation}
\label{limp}
\Big(\frac{Y^{N,M}(t)-Mt\ide_N}{\sqrt{M}}\Big)_{t\ge 0} \overset{d}{\underset{M \to
    \infty}{\longrightarrow}} (X^N(t^2))_{t\ge 0}
\end{equation}
in the sense of weak convergence in $\cC(\dR_+,\cH_N)$.
\end{thm}

As announced, the proofs of the previous theorems are postponed up to section
(\ref{proofs}). Their ideas are quite simple: for Theorem \ref{lim},
we can compute directly on the densities given by (\ref{densg}) and
(\ref{densl}). For Theorem \ref{processus}, our central-limit
convergence is shown to follow from a law of large numbers at the
level of quadratic variations. We use the usual two-stepped argument
for the convergence of processes: convergence of finite-dimensionnal distributions and tightness.

Let us mention the straightforward consequence of Theorems \ref{lim} and
\ref{processus} on the convergence of eigenvalues. If $H \in \cH_N$,
let us denote by $l_1(H) \le \cdots \le l_N(H)$ its (real) eigenvalues
and $l(H)=(l_1(H),\ldots,l_N(H))$. Using the min-max formulas, it is
not difficult to see that each $l_i$ is $1$-Lipschitz for the
Euclidean norm on $\cH_N$. Thus, $l$ is continuous on $\cH_N$. Therefore, if we
set $\mu^{N,M}=l(Y^{N,M})$ and $\la^N=l(X^N)$ 
\begin{equation}
\label{vps}
\Big(\frac{\mu_i^{N,M}-M}{\sqrt{M}}\Big)_{1\le i \le N} \overset{d}{\underset{M \to
    \infty}{\longrightarrow}} (\la_i^N)_{1\le i \le N} 
\end{equation}
With the obvious notations, the process version also takes place:
\begin{equation}
\label{vpp}
\Big( \Big(\frac{\mu_i^{N,M}(t)-Mt}{\sqrt{M}}\Big)_{1\le i \le N}\Big)_{t\ge 0} \overset{d}{\underset{M \to
    \infty}{\longrightarrow}} \big((\la_i^N(t^2))_{1\le i \le N}\big)_{t\ge 0}
\end{equation}

Analogous results hold in the real case of GOE and LOE and they can be proved with the same arguments.

For connections with the results of this section, see Theorem 2.5 of
\cite{dette} and a note in Section 5 of \cite{oconnell-brownian}. The basic form of our central-limit theorem already appeared without proof in the
Introduction of \cite{jonsson}. To our knowledge, the process version had not
been considered in the existing literature.

\section{Consequences on representations for eigenvalues}
\label{int}

\subsection{The largest eigenvalue}
\label{maxim}

Let us first indicate how to recover from (\ref{vps}) the identity 
\begin{equation}
\label{iden}
\lam \overset{d}{=} \sup_{0=t_0\le\ldots\le t_N =1} \sum\limits_{i=1}^{N}
(B_i(t_i)-B_i(t_{i-1}))
\end{equation} 
where $\lam=\la_N^N$ is the maximal eigenvalue of $\gue$ and
$(B_i\,,\,1\le i\le N)$ is a standard $N$-dimensional Brownian
motion. If $(w_{i,j}\,,\,(i,j) \in (\dN \setminus \{0\})^2)$ are i.i.d. exponential variables with parameter one, define 
\begin{equation}
\label{defi}
H(M,N)=\max \Big\{ \sum\limits_{(i,j)\in \pi}
w_{i,j}\; ;\; \pi \in \cP(M,N) \Big\}
\end{equation}
where $\cP(M,N)$ is the set of all paths $\pi$ taking
only unit steps in the north-east direction in the rectangle
$\{1,\ldots,M\}\times \{1,\ldots,N\}$. In \cite{johansson-shape}, it is
noticed that 
\begin{equation}
\label{rema}
H(M,N) \,\overset{d}{=} \, \mu_{\text{max}}^{M,N}
\end{equation}
where $\mu_{\text{max}}^{M,N}=\mu_N^{N,M}$ is the largest eigenvalue of
$\lue$. Now an invariance principle due to Glynn and Whitt in \cite{gly-whi} shows that 
\begin{equation}
\label{gl}
\frac{H(M,N)-M}{\sqrt{M}} \overset{d}{\underset{M
    \to \infty}{\longrightarrow}} \sup_{0=t_0\le\ldots\le t_N =1}
\sum\limits_{i=1}^{N} (B_i(t_i)-B_i(t_{i-1})) \; .
\end{equation}
On the other hand, by (\ref{vps}) 
\begin{equation}
\label{vpmax}
\frac{\mu_{\text{max}}^{N,M}-M}{\sqrt{M}} \overset{d}{\underset{M \to
    \infty}{\longrightarrow}} \lam \; .
\end{equation}
Comparing (\ref{rema}), (\ref{gl}) and (\ref{vpmax}), we get
(\ref{iden}) for free.

In the next section, we will give proofs of more general statements than (\ref{rema}) and (\ref{gl}).

\subsection{The other eigenvalues}
\label{other}
In fact, Johansson's observation involves all the eigenvalues of
$\lue$ and not only the largest one. Although it does not appear
exactly like that in \cite{johansson-shape}, it takes the following form. First, we
need to extend definition (\ref{defi}) as follows: for each $k$, $1\le k \le N$, set
\begin{equation}
\label{toutesvp} 
H_k(M,N)=\max \Big\{\sum\limits_{(i,j)\in \pi_1 \cup \cdots \cup \pi_k}
w_{i,j}\; ;\; \pi_1,\ldots ,\pi_k \in \cP(M,N)\, ,\,
\pi_1,\ldots,\pi_k \,\text{all disjoint} \;\Big\}\; .
\end{equation}
Then, the link, analogous to (\ref{rema}), with the eigenvalues of $\lue$ is
expressed by  
\begin{equation}
\label{toutesvp2}
H_k(M,N)\,\overset{d}{=} \,
\mu_N^{N,M}+\mu_{N-1}^{N,M}+\cdots+\mu_{N-k+1}^{N,M} \; .
\end{equation}
In fact, the previous equality in distribution is also valid for the vector $(H_k(M,N))_{1\le k \le N}$ and the corresponding sums of
eigenvalues, which gives a representation for all the eigenvalues of
$\lue$. 

\begin{proof}[Proof of (\ref{toutesvp2})]
The arguments and notations are taken from Section $2.1$ in \cite{johansson-shape}. Denote
by $\cM_{M,N}$ the set of $M \times N$ matrices $A=(a_{ij})$ with non-negative
integer entries and by $\cM_{M,N}^s$ the subset of $A \in \cM_{M,N}$ such that
$\Si(A)=\sum a_{ij}=s$. Let us recall that the Robinson-Schensted-Knuth
(RSK) correspondence is a one-to-one mapping from $\cM_{M,N}^s$ to the set
of pairs $(P,Q)$ of semi-standard Young tableaux of the same shape $\la$ which is a
partition of $s$, where $P$ has elements in $\{1,\ldots ,N\}$ and $Q$ has
elements in $\{1,\ldots ,M\}$. Since $M \ge N$ and since the numbers are strictly
increasing down the columns of $P$, the number of rows of $\la$ is at most
$N$. We will denote by $\text{RSK}(A)$ the pair of
Young tableaux associated to a matrix $A$ by the RSK correspondence and by
$\la(\text{RSK}(A))$ their commun shape. The crucial fact about this correspondence is the combinatorial property that, if
$\la=\la(\text{RSK}(A))$, then for all $k$, $1 \le k \le N$,
\begin{equation}
\label{aqwx}
\la_1 + \la_2 +\cdots +\la_k = \max \Big\{\sum\limits_{(i,j)\in \pi_1 \cup \cdots \cup \pi_k} a_{i,j}\; ;\; \pi_1,\ldots ,\pi_k \in \cP(M,N)\, ,\,
\pi_1,\ldots,\pi_k \,\text{all disjoint} \; \Big\}\, .
\end{equation}
Now consider a random $M \times N$ matrix $X$ whose entries $(x_{ij})$ are i.i.d. geometric variables with parameter $q$. Then for any $\la^0$ partition
of an integer $s$, we have 
$$
\dP(\la(\text{RSK}(X))=\la^0\,) = \sum\limits_{A \in \cM_{M,N}^s
  \; ,\;\la(\text{RSK}(A))=\la^0} \dP(X=A) \, .$$
But for $A \in \cM_{M,N}^s$, $\dP(X=A)=(1-q)^{MN}q^s$ is independent of $A$, 
which implies $$\dP(\la(\text{RSK}(X)=\la^0))=(1-q)^{MN}q^{\sum
  \la_i^0}\,L(\la^0,M,N)$$ where $L(\la^0,M,N)=\sharp \{ A \in \cM_{M,N}\,
,\,\la(\text{RSK}(A))=\la^0 \}$. Since the RSK mapping is one-to-one
$$L(\la^0,M,N)=Y(\la^0,M)\,Y(\la^0,N)$$ where $Y(\la^0,K)$ is just the number
of semi-standard Young tableaux of shape $\la^0$ with elements in $\{1,\ldots
,K\}$. This cardinal is well-known in combinatorics and finally 
$$L(\la^0,M,N)=c_{MN}^{-1} \, \prod\limits_{1 \le i < j \le N} (h_j^0-h_i^0)^2 \,
\prod\limits_{1 \le i \le N} \frac{(h_i^0+M-N)!}{h_i^0!}$$
where $c_{MN}=\prod\limits_{0 \le i \le N-1} j\,! \, (M-N+j)\,!$ and $h_i^0=\la_i^0+N-i$ such that $h_1^0 > h_2^0 > \cdots > h_N^0 \ge 0$. With the same correspondence as before
between $h$ and $\la$, we can write
\begin{eqnarray*}
\dP(h(\text{RSK}(X))=h^0) &=&c_{MN}^{-1}
\,\frac{(1-q)^{MN}}{q^{N(N-1)/2}}\,\prod\limits_{1 \le i < j \le N} (h_j^0-h_i^0)^2
\,\prod\limits_{1 \le i \le N} \frac{(h_i^0+M-N)!}{h_i^0!} \\
&\overset{\text{def}}{=}& \rho_{(M,N,q)}(h^0)\, .
\end{eqnarray*}
Now set $q=1-L^{-1}$ and use the notation $X_L$ instead of $X$ to recall the
dependence of the distribution on $L$. An easy asymptotic expansion shows that $$
L^{N} \rho_{(M,N,1-L^{-1})}(\lfloor Lx \rfloor ) \underset{L \to
  \infty}{\longrightarrow} d_{MN}^{-1} \, \prod\limits_{1 \le i < j \le N}
(x_j-x_i)^2 \, \prod\limits_{1 \le i  \le N}  x_i^{M-N} e^{-x_i}=\rho_{\lue}(x)$$ 
where $\rho_{\lue}$ is the joint density of the ordered eigenvalues of $\lue$. This can be used to prove that 
\begin{equation}
\label{azer}
\frac{1}{L} h(\text{RSK}(X_L)) \overset{d}{\underset{L \to
  \infty}{\longrightarrow}} (\mu_N^{MN},\mu_{N-1}^{MN},\ldots,\mu_1^{MN})\, . 
\end{equation}
On the other hand, if $x_L$ is a geometric variable with parameter $1-L^{-1}$,
then $x_L /L$ converges in ditribution, when $L \to \infty$, to an exponential
variable of parameter one. Therefore, using the link between $h$ and $\la$
together with (\ref{aqwx}), we have
$$\frac{1}{L} \Big(\sum\limits_{i=1}^k h_i(\text{RSK}(X_L)) \Big)_{1 \le k\le
  N} \overset{d}{\underset{L \to  \infty}{\longrightarrow}} (H_k(M,N))_{1 \le
  k\le N}\, .$$
Comparing with (\ref{azer}), we get the result.
\end{proof}

Now, let us try to adapt what we previously did with $H(M,N)$
and $\mum$ to the new quantities $H_k(M,N)$. First, we would like to have an analogue of the Glynn-Whitt invariance principle (\ref{gl}). To
avoid cumbersome notations, let us first look at the case
$k=2,\,N=3$. In this case, the geometry involved in the $H_2(M,3)$ is
simple: we are trying to pick up the largest possible weight by using
two north-east disjoint paths in the rectangle $\{1,\ldots,M\}\times
\{1,2,3\}$. The most favourable configuration corresponds to one
path (the bottom one) starting at $(1,1)$ and first going right. Then
it jumps to some point of $\{2,\ldots,M\}\times \{2\}$ and goes horizontally up to
$(M,2)$. The upper path starts at $(1,2)$, will also jump and go right
up to $(M,3)$. The constraint that our paths must be disjoint forces the
$x$-coordinate of the jump of the bottom path to be larger than that
of the jump of the upper path. This corresponds to the obvious figure
1.

\begin{figure}[htbp]
\begin{center}
\input{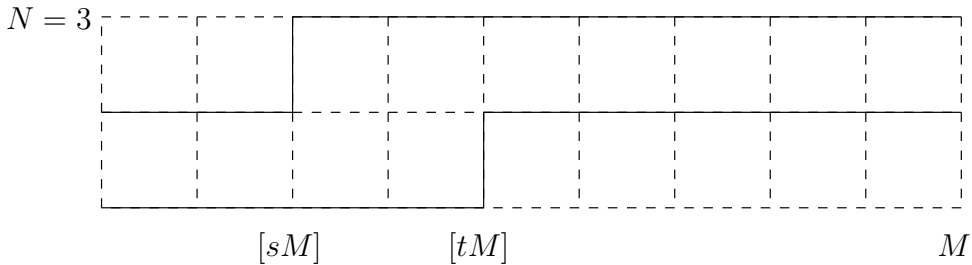}
\caption{Configuration of paths in the case $k=2$ and $N=3$}
\end{center}
\end{figure}

This figure suggests that in the Donsker limit of random walks
converging to Brownian motion, we will have  
$$\frac{H_2(M,3)-2M}{\sqrt{M}} \overset{d}{\underset{M
    \to \infty}{\longrightarrow}} \Om_2^{(3)} \overset{\text{def}}{=} \sup_{0\le s \le t\le 1} (\,B_1(t)+B_2(s)+B_2(1)-B_2(t)+B_3(1)-B_3(s)\,)$$
where $(B_1,B_2,B_3)$ is standard $3$-dimensional Brownian motion.

For the case of $k=2$ and general $N$, we have the same configuration
except that the number of jumps for each path will be $N-2$ so that 
\begin{equation}
\label{edc}
\frac{H_2(M,N)-2M}{\sqrt{M}} \overset{d}{\underset{M
    \to \infty}{\longrightarrow}} \Om_2^{(N)} \overset{\text{def}}{=} \sup \,\sum\limits_{j=1}^N \,(\,B_j(s_{j-1})-B_j(s_{j-2})+B_j(t_j)-B_j(t_{j-1})\,)
\end{equation}
where $(B_j)_{1\le j \le N}$ is a standard $N$-dimensional Brownian
motion and the $\sup$ is taken over all subdivisions of $[0,1]$ of the
following form: $$ 0=s_{-1}=s_0=t_0 \le s_1\le t_1 \le s_2\le t_2 \le \cdots \le  t_{N-2}\le s_{N-1}=t_{N-1}=s_N=t_N=1 \; . $$

\begin{proof}[Proof of limit (\ref{edc})]

Let us first consider the case of $H_2(M,N)$ :
$$H_2(M,N)=\max \{\sum\limits_{(i,j)\in \pi_1 \cup \pi_2}
w_{i,j}\; ;\; \pi_1,\,\pi_2  \in \cP(M,N)\; ;\; \pi_1,\,\pi_2 \; \text{disjoint} \}$$
Since our paths are disjoint, one (say $\pi_1$) is always lower than
the other (say $\pi_2$): for all $i \in \{1,\ldots,M\}$,  $\max \{j\, ;\, (i,j) \in
\pi_1\} < \min \{j\, ;\, (i,j) \in \pi_2 \}$. We will denote this by
$\pi_1 < \pi_2$. Then, it is not difficult to see on a picture that,
for any two paths $\pi_1 < \pi_2 \, \in \cP(M,N)$, one can always
find paths $\pi_1'< \pi_2' \, \in \cP(M,N)$ such that $\pi_1 \cup
\pi_2 \subset \pi_1' \cup \pi_2'$ , $\pi_1'$ starts from $(1,1)$, visits $(2,1)$ then finishes in $(M,N-1)$ and $\pi_2'$ starts from
$(2,1)$ and goes up to $(M,N)$. Let us call $\cP(M,N)'$ the set of
pairs of such paths $(\pi_1',\pi_2')$. Thus 
$$H_2(M,N)=\max \{\sum\limits_{(i,j)\in \pi_1 \cup \pi_2}
w_{i,j}\; ;\; (\pi_1,\pi_2) \in \cP(M,N)'\,\} \, .$$
Now two paths $(\pi_1,\pi_2) \in \cP(M,N)'$ are uniquely determined by
the non-decreasing sequences of their $N-2$ vertical jumps, namely
$0\le t_1 \le \ldots \le t_{N-2}\le 1$ for $\pi_1$ and $0\le s_1 \le
\ldots \le s_{N-2}\le 1$ for $\pi_2$ such that:
\begin{itemize}
\item[-] $\pi_1$ is horizontal on $\left[\, \lfloor t_{i-1}M \rfloor
    ,\lfloor t_i M \rfloor \, \right]
  \times \{i\}$ and vertical on $\{\lfloor t_i M\rfloor \}\times [i,i+1]$,
\item[-] $\pi_2$ is horizontal on $\left[\, \lfloor s_{i-1}M\rfloor
    ,\lfloor s_i M\rfloor \, \right]
  \times \{i+1\}$ and vertical on $\{\lfloor s_i M\rfloor \}\times [i+1,i+2]$,
\item[-] $s_i < t_i$ for all $i \in \{1,\ldots,N-2\}$, this constraint being
  equivalent to the fact that $\pi_1 < \pi_2$ .
\end{itemize}
The weight picked up by two such paths coded by $(t_i)$ and $(s_i)$ is
\begin{itemize}
\item[-] $w_{1,1}+w_{2,1}+\cdots+w_{\lfloor t_1 M \rfloor ,1}$ on the
  first floor,
\item[-] $w_{1,2}+\cdots+w_{\lfloor s_1 M \rfloor ,2}\; + \;w_{\lfloor
    t_1 M \rfloor ,2}+\cdots+w_{\lfloor t_2 M \rfloor ,2} $ on the second
  floor,
\item[-] $w_{\lfloor s_1 M \rfloor ,3}+\cdots+w_{\lfloor s_2 M \rfloor
    ,3}\; + \; w_{\lfloor
    t_2 M \rfloor ,3}+\cdots+w_{\lfloor t_3 M \rfloor ,3} $ on the
  third floor,
\item[-] and so on, up to floor $N$ for which the contribution is
  $w_{\lfloor s_{N-2} M \rfloor ,N}+\cdots+w_{M,N}$ .
\end{itemize}

This yields 
$$H_2(M,N)=\sup \; \sum\limits_{j=1}^N \,\bigg( \,\sum\limits_{i=\lfloor s_{j-2}M \rfloor}^{\lfloor s_{j-1}M
  \rfloor} w_{i,j} \; + \; \sum\limits_{i=\lfloor t_{j-1}M \rfloor}^{\lfloor
  t_{j}M \rfloor} w_{i,j}\;\bigg)\, .$$ 
Hence,
\begin{equation*}
\begin{split}
\frac{H_2(M,N)-2M}{\sqrt{M}} =\sup \; \sum\limits_{j=1}^N &\,\bigg( \,\frac{ \sum\limits_{i=\lfloor s_{j-2}M \rfloor}^{\lfloor s_{j-1}M
  \rfloor} w_{i,j} -(s_{j-1}-s_{j-2})M}{\sqrt{M}} \\
&\; + \;\frac{\sum\limits_{i=\lfloor t_{j-1}M \rfloor}^{\lfloor
  t_{j}M \rfloor} w_{i,j} -(t_j-t_{j-1})M}{\sqrt{M}}\;\bigg)\, .
\end{split}
\end{equation*}
Donsker's principle states that 
$$\bigg( \frac{\sum\limits_{i=1}^{\lfloor  sM \rfloor} w_{ij} - sM}{\sqrt{M}}
\bigg)_{1\le j \le N} \overset{d}{\underset{M \to  \infty}{\longrightarrow}} (B_j(s))_{1 \le j\le N}$$
where the convergence takes place in the space of
cadlag trajectories of the variable $s \in\dR_+$ equipped with the Skorohod
topology. This allows us to conclude (see \cite{gly-whi} for a detailed account on the continuity of our mappings in the Skorohod topology).
\end{proof}

For general $k$ and $N$, the same pattern works with $k$ disjoint paths having each $N-k$
jumps. This yields the following
central-limit behaviour: 
\begin{equation}
\label{lim4}
\frac{H_k(M,N)-kM}{\sqrt{M}} \overset{d}{\underset{M
    \to \infty}{\longrightarrow}} \Om_k^{(N)} \overset{\text{def}}{=} \sup
\,\sum\limits_{j=1}^N \, \sum\limits_{p=1}^k (\,B_j(s_{j-p+1}^p)-B_j(s_{j-p}^p)\,)
\end{equation}
where the $\sup$ is taken over all subdivisions $(s_i^p)$ of $[0,1]$ of the
following form: $$s_i^p\in [0,1]\; ,\; s_i^{p+1} \le s_i^p \le s_{i+1}^p \; ,\; s_i^p=0 \;\text{for}\; i\le 0\;
\text{and}\; s_i^p=1\;\text{for}\; i\ge N-k+1 $$

Now, imitating the argument for the $\lam$, we obtain that
\begin{equation}
\label{perso}
 \Om_k^{(N)} \overset{d}{=} \la_N^N + \la_{N-1}^N + \cdots +\la_{N-k+1}^N
\end{equation}
where we recall that $\la_1^N \le \cdots \le \la_N^N$ are the
eigenvalues of $\gue$. In fact, the previous equality is also true
when considering the vector $(\Om_k^{(N)})_{1\le k\le N}$ and the
corresponding sums of eigenvalues, which yields a representation for
all the eigenvalues of $\gue$.\\

A representation for the eigenvalues of $\gue$ was already
obtained in \cite{oconnell-representation}. Let us compare both representations. Denote by $\cD_0(\dR_+)$ the space of cadlag paths $f : \dR_+ \to \dR$ with
$f(0)=0$ and for $f,g \in \cD_0(\dR_+)$, define $f  \ot g \in
\cD_0(\dR_+)$ and $f\od g \in \cD_0(\dR_+)$ by 
$$f\ot g(t)=\,\inf\limits_{0\le s\le t} \left( f(s)+g(t)-g(s) \right) \quad
\text{and} \quad f\od g(t)=\,\sup\limits_{0\le s\le t} \left( f(s)+g(t)-g(s) \right)$$
By induction on $N$, define $\Ga^{(N)} : \cD_0(\dR_+)^N \to
\cD_0(\dR_+)^N$ by $$\Ga^{(2)}(f,g)=(f\ot g,g\od f)$$ and for $N>2$ and
$f=(f_1,\ldots,f_N)$ $$\Ga^{(N)}(f)=\left(f_1 \ot \cdots \ot f_N,
\Ga^{(N-1)}(f_2 \od f_1,f_3\od (f_1\ot f_2),\ldots,f_N\od (f_1\ot \cdots
\ot f_{N-1}))\right) \, .$$
Then the main result in \cite{oconnell-representation} is : 
\begin{equation}
\label{yor}
\la^N \overset{d}{=} \Ga^{(N)}(B)(1)
\end{equation}
where $B=(B_i)_{1\le i \le N}$ is standard $N$-dimensional Brownian
motion and $\la^N$ is the vector of eigenvalues of $\gue$. 
In fact, it is proved in \cite{oconnell-representation} that
identity (\ref{yor}) is true for the whole processes and not only
their marginals at time $1$.

Thus $$\la_N^N + \la_{N-1}^N + \cdots +\la_{N-k+1}^N \overset{d}{=}
\Ga^{(N)}_N(B)(1) + \Ga^{(N)}_{N-1}(B)(1) +\cdots + \Ga^{(N)}_{N-k+1}(B)(1)\, .$$ 
Comparison with (\ref{perso}) gives   
\begin{equation}
\label{bizarre}
\Om_k^{(N)} \overset{d}{=}
\Ga^{(N)}_N(B)(1) + \Ga^{(N)}_{N-1}(B)(1) +\cdots + \Ga^{(N)}_{N-k+1}(B)(1) \,
.
\end{equation}
This equality in distribution also holds for the $N$-vector $(\Om_k^{(N)})_{1 \le k \le N}$.

Now let us remark that the definition of the components $\Ga^{(N)}_k$ of
$\Ga^{(N)}$ is quite intricate: it involves a sequence of nested ``$\inf$'' and
``$\sup$''. On the contrary, $\Om_k^{(N)}$ is only defined by one ``$\sup$''
but over a complicated sequence of nested subdivisions. We ignore whether
these identities are: trivial and uninteresting; already well-known;
 true for the deterministic formulas (ie true when replacing
  independent Brownian motions by continuous functions) or true only
  in distribution.

Our concern raises the question about the link between the $\Ga^{(N)}$
introduced in \cite{oconnell-representation} and the
Robinson-Schensted-Knuth correspondence that gave birth to our
$\Om^{(N)}$. A striking equivalence between both objects was recently
put forward in the context of random words (\cite{oconnell-path-transformation}).  

Finally, let us notice that the heart of our arguments to get the
previous representations is the identity (\ref{toutesvp}). The proof presented
here is taken from \cite{johansson-shape} and is organized in two
steps : first the computation of the joint density for $(H_k(M,N))_{1
  \le k \le N}$ by combinatorial means and second the observation that
this density coincides with the eigenvalue density of $\lue$. It would
be tempting to get a deeper understanding of this result. This would
all amount to obtaining a representation for non-colliding squared Bessel processes.

\section{Proofs}
\label{proofs}

\begin{proof}[Proof of Theorem \ref{lim}]

The density of $Y^{N,M}$ on the space $\cH_N$ is given by $\ffi_M(S) \II_{S \ge 0}dS$ where $$\ffi_M(S)=\frac{1}{\pi^{N(N-1)/2} \, \prod_{j=1}^N
  (M-j)!}\, (\det S)^{M-N} \exp (-\Tr S) \, .$$
We make the following change of variables :
$H=\frac{S-M\ide_N}{\sqrt{M}}$. The density $\psi_M$ of
$\frac{Y^{N,M}-M\ide_N}{\sqrt{M}}$ is obtained from
$\psi_M(H)dH=\ffi_M(S)dS$ and can be written  
\begin{eqnarray*}
\psi_M(H) & = & M^{N^2/2}\ffi_M(M\ide_N+\sqrt{M}H) \\
& = & \frac{M^{N^2/2} M^{N(M-N)}}{\pi^{N(N-1)/2} \, \prod_{j=1}^N
  (M-j)!}\, \det \Big(\ide_N+\frac{H}{\sqrt{M}}\Big)^{M-N}\,
\exp\Big(-M\Tr\big(\ide_N+\frac{H}{\sqrt{M}}\big)\Big)\\
& = & c_M \exp \Big\{ (M-N)\log \det (\ide_N+\frac{H}{\sqrt{M}}) -M\Tr
\frac{H}{\sqrt{M}} \Big) 
\end{eqnarray*}
where $c_M=\exp(-MN) \frac{M^{N^2/2} M^{N(M-N)}}{\pi^{N(N-1)/2} \,
  \prod_{j=1}^N (M-j)!}$. But if $(\la_i)$ are the eigenvalues of $H
\in \cH_N$, then for small $\eps$ : 
\begin{eqnarray*}
\det(\ide_N+\eps H) & = & 1+\eps \sum_i \la_i +\eps^2 \sum_{i<j} \la_i
\la_j +\, \cO(\eps^3)\\
& = & 1+\eps \Tr H +\frac{\eps^2}{2}[(\Tr H)^2 -\Tr H^2] + \cO(\eps^3)
\end{eqnarray*}
 Thus $$\log \det(\ide_N+\eps H) = \eps \Tr H -\frac{\eps^2}{2} \Tr H^2 +\cO(\eps^3)$$
which results in 
\begin{equation}
\label{asd}
\psi_M(H)= c_M \exp \{ (M-N)(\frac{\Tr
  H}{\sqrt{M}}-\frac{1}{2M}\Tr H^2 + \cO(\frac{1}{M^{3/2}}))-M \Tr
\frac{H}{\sqrt{M}} \}
\end{equation}
The exponential term in (\ref{asd}) converges to $\exp(-\frac{1}{2}\Tr H^2)$ as $M
\to \infty$. Then
using Stirling's formula and the fact that $\prod_{j=1}^N (M-j)! \underset{M \to
  \infty}{\sim} \frac{(M!)^N}{M^{N(N+1)/2}}$, we
find that $c_M \underset{M \to  \infty}{\to} \frac{1}{2^{N/2}
  \pi^{N^2/2}}$. Thus the density $\psi_M$ converges to the density
$\psi$ of $\gue$. In fact, it is not difficult to see that one can have sufficient control of the rest in the previous asymptotic expansions so as to obtain uniform boundedness of $\psi_M$ on compact subsets of $\cH_N$. Therefore, dominated convergence yields $$ \int\limits_{\cH_N}
f(H)\psi_M(H)\,dH \underset{M \to \infty}{\to} \int\limits_{\cH_N}
f(H)\psi(H)\,dH$$
for every $f$ continuous, compactly supported on $\cH_N$.
\end{proof}

\begin{proof}[Proof of Theorem \ref{processus}]

We will write $A$ instead of $A^{M,N}$. For $1 \le i\le N \, , \, 1\le j\le M$, the superscript $ij$ when applied to a matrix stands for its entry at line
$i$ and column $j$. The value at time $t$ of any process $x$ will be denoted either
$x(t)$ or $x_t$. Let us set
$$Z_M(t)=\frac{Y^{N,M}(t)-Mt\ide_N}{\sqrt{M}}=\frac{AA^*(t)-Mt\ide_N}{\sqrt{M}}\,
.$$
Then $$Z_M^{ij}=\frac{1}{\sqrt{M}} (\sum\limits_{k=1}^M A^{ik} \overline{A}^{jk}
-Mt\de_{ij}) \quad , \quad dZ_M^{ij}=\frac{1}{\sqrt{M}}
\sum\limits_{k=1}^M (A^{ik} d\overline{A}^{jk} + \overline{A}^{jk} dA^{ik})\, ,$$
which implies $$dZ_M^{ij}\cdot dZ_M^{i'j'} =\frac{1}{M}
\sum\limits_{k=1}^M (A^{ik} \overline{A}^{j'k}\de_{i'j} + \overline{A}^{jk}
A^{i'k} \de_{ij'})\, dt\, .$$
The quadratic variation follows to be : $$\langle Z_M^{ij},Z_M^{i'j'} \rangle_t=\frac{1}{M}
\sum\limits_{k=1}^M \int_0^t (A^{ik}_s \overline{A}^{j'k}_s\de_{i'j} + \overline{A}^{jk}_s
A^{i'k}_s \de_{ij'})\, ds\, .$$
By the classical law of large numbers, we get that this converges
almost surely to :
$$\int_0^t \left(\dE(A^{i1}_s \overline{A}^{j'1}_s)\de_{i'j} + \dE(\overline{A}^{j1}_s
A^{i'1}_s) \de_{ij'}\right)\, ds\;=\;\int_0^t \de_{ij'}\de_{i'j}
2s\,ds\;=\;t^2\de_{ij'}\de_{i'j}\, .$$

Note that he previous formula shows that, in the limit, the quadratic variation is $0$
if $i\ne j'$ and $i' \ne j$, which is obvious even for finite $M$ without
calculations. However, if for instance $i=j'$ and $i' \ne j$, then the
quadratic variation is not $0$ for finite $M$ and only becomes null in
the limit. This is some form of asymptotic independence.

First, let us prove tightness of the process $Z_M$ on any fixed
finite interval of time $[0,T]$. It is sufficient to prove tightness
for every component, let us do so for $Z_M^{11}$ for example
($Z_M^{11}$ is real). We will apply Aldous' criterion (see \cite{kipnis-landim}). Since
$Z_M^{11}(0)=0$ for all $M$, it is enough to check that, for all $ \eps >0$, 
\begin{equation}
\label{ald}
\lim\limits_{\de \to 0}\, \lsup\limits_{M \to \infty} \,
\sup\limits_{\tau \,,\, 0\le \te \le \de} \dP(\, |Z_M^{11}(\tau
+\te)-Z_M^{11}(\tau)|\ge \eps\,)\,=\,0
\end{equation}
where the $\sup$ is taken over all stopping times $\tau$ bounded by
$T$. For $\tau$ such a stopping time, $\eps > 0$ and $0\le \te \le \de \le 1$, we have
\begin{eqnarray*}
\dP(\, |Z_M^{11}(\tau +\te)-Z_M^{11}(\tau)|\ge \eps\,) & \le &
\frac{1}{\eps^2}\dE ((Z_M^{11}(\tau +\te)-Z_M^{11}(\tau))^2)\\
& = & \frac{1}{\eps^2}\dE (\int\limits_{\tau}^{\tau +\te} d\langle Z_M^{11},Z_M^{11} \rangle_t)\\
& = & \frac{2}{M\eps^2} \sum\limits_{k=1}^M \dE
(\int\limits_{\tau}^{\tau +\te} |A^{1k}_s|^2\,ds)\\
& \le & \frac{2}{M\eps^2} \sum\limits_{k=1}^M
\dE(\te \,\sup\limits_{0\le s \le T+1} |A^{1k}_s|^2)\\
& = & \frac{2\te}{\eps^2}  \dE(\sup\limits_{0\le s \le T+1}
|A^{11}_s|^2)
\end{eqnarray*}
Since $c_T= \dE(\sup\limits_{0\le s \le T+1}\,|A^{11}_s|^2 ) < \infty$, then 
$$\lsup\limits_{M \to \infty} \,
\sup\limits_{\tau \,,\, 0\le \te \le \de} \dP(\, |Z_M^{11}(\tau
+\te)-Z_M^{11}(\tau)|\ge \eps\,)\,\le\,\frac{2\de\, c_T}{\eps^2}\, .$$
This last line obviously proves (\ref{ald}).

Let us now see that the finite-dimensionnal distributions converge to the appropriate limit. Let us first fix $i,j$ and look at the component
$Z_M^{ij}=\frac{x_M+\sqrt{-1}y_M}{\sqrt{2}}$. We can write 
\begin{equation}
\label{formula}
\langle x_M,y_M \rangle_t =0 \quad , \quad
\langle x_M,x_M \rangle_t = \langle y_M,y_M \rangle_t = \frac{1}{M}
\sum\limits_{k=1}^M \int_0^t \al_s^k \,ds
\end{equation}
where $\al_s^k=|A^{ik}_s|^2 +|A^{jk}_s|^2$. We are going to consider
$x_M$. Let us fix $T \ge 0$. For any $(\nu_1,\ldots ,\nu_n) \in [-T,T]^n$
and any $0=t_0 < t_1
< \ldots <t_n \le T$,
we have to prove that  
\begin{equation}
\label{qaz}
\dE \bigg( \exp \Big( i \sum\limits_{j=1}^n \nu_j  (x_M(t_j)-x_M(t_{j-1})\Big)\bigg) \underset{M \to
    \infty}{\longrightarrow} \exp \Big( \sum\limits_{j=1}^n \frac{\nu_j^2}{2}
  (t_{j}^2-t_{j-1}^2)\Big)\, .
\end{equation}
We can
always suppose $|t_j-t_{j-1}| \le \de$ where $\de$ will be chosen
later and will only depend on $T$ (and not on $n$). We will prove property (\ref{qaz}) by induction on $n$. For $n=0$, there is nothing to
prove. Suppose it is true for $n-1$. Denote by $(\cF_t)_{t\ge 0}$ the
filtration associated to the process $A$. Then write: 
\begin{equation}
\label{hypoa}
\dE \Big( e^{i \sum\limits_{j=1}^n \nu_j
(x_M(t_j)-x_M(t_{j-1}))}\Big)=\dE\Big(e^{i \sum\limits_{j=1}^{n-1}
\nu_j (x_M(t_j)-x_M(t_{j-1}))}\,\dE\big(e^{i(x_M(t_n)-x_M(t_{n-1}))} |
\cF_{t_{n-1}}\big)\Big)\, .
\end{equation}
We define the martingale $\cM_t=e^{i\nu_n x_M(t)-\frac{\nu_n^2}{2} \langle
  x_M,x_M \rangle_t}$. Hence $$\dE \big(e^{i\nu_n (x_M(t_n)-x_M(t_{n-1}))} |
\cF_{t_{n-1}}\big)=\dE \left(\frac{\cM_{t_n}}{\cM_{t_{n-1}}}
e^{\frac{\nu_n^2}{2} \langle
  x_M,x_M \rangle_{t_{n-1}}^{t_n}}\, |\,\cF_{t_{n-1}}\right)$$ 
with the notation $\langle  x,x \rangle_s^t =\langle  x,x \rangle_t - \langle  x,x \rangle_s $. This yields 
\begin{equation}
\label{qwer}
 e^{-\frac{\nu_n^2}{2}(t_n^2-t_{n-1}^2)} \dE\big(e^{i\nu_n(x_M(t_n)-x_M(t_{n-1}))}
\, |\, \cF_{t_{n-1}}\big) -1 
=\dE
\left( \frac{\cM_{t_n}}{\cM_{t_{n-1}}}\, \zeta_M\;
  |\,\cF_{t_{n-1}}  \right) 
\end{equation}
where we set $\zeta_M=e^{\frac{\nu_n^2}{2}( \langle x_M,x_M
  \rangle_{t_{n-1}}^{t_n}-(t_n^2-t_{n-1}^2))}-1$. Using that $|e^z-1|\le |z| e^{|z|}$, we deduce that
$$|\zeta_M| \le K\, | \langle x_M,x_M
  \rangle_{t_{n-1}}^{t_n}-(t_n^2-t_{n-1}^2)| \, e^{\frac{\nu_n^2}{2} \langle x_M,x_M
  \rangle_{t_{n-1}}^{t_n}}$$ where $K=\nu_n^2 /2$. The Cauchy-Schwarz
inequality implies that
$$\dE (|\zeta_M|) \le K\,
\left( \dE \left( \langle x_M,x_M
  \rangle_{t_{n-1}}^{t_n}-(t_n^2-t_{n-1}^2)\right)^2 \right)^{1/2}\,
\left( \dE\left( e^{\nu_n^2 \langle x_M,x_M  \rangle_{t_{n-1}}^{t_n}}\right)
\right)^{1/2}\, .$$ 
By convexity of the function $x \to e^x$ : $$e^{\nu_n^2 \langle x_M,x_M \rangle
    _{t_{n-1}}^{t_n}}\, =\, \exp\left(\frac{1}{M}
\sum\limits_{k=1}^M \nu_n^2 \int_{t_{n-1}}^{t_n} \al_u^k \,du\right)\,\le\,\frac{1}{M}
\sum\limits_{k=1}^M e^{\nu_n^2\,(t_n-t_{n-1})\sup\limits_{0\le u\le
    t_n} \al_u^k}$$
and thus
$$\dE \left(e^{\nu_n^2 \langle x_M,x_M \rangle
    _{t_{n-1}}^{t_n}}\right)\, \le \,\frac{1}{M}
\sum\limits_{k=1}^M \dE\left(e^{\nu_n^2\,(t_n-t_{n-1})\sup\limits_{0\le u\le
    t_n} \al_u^k}\right)\, =\,\dE\left(e^{\nu_n^2(t_n-t_{n-1})
  \sup\limits_{0\le u\le t_n} \al_u^1}\right)\, .$$
Now let us recall that $\al_u^1=|A^{i1}_u|^2 +|A^{j1}_u|^2$, which
means that $\al^1$ has the same law as a sum of squares of four
independent Brownian motions. It is then easy to see that there exists
$\de>0$ (depending only on $T$) such that $\dE(\exp\, (T^2 \de
  \sup\limits_{0\le u\le T} \al_u^1))< \infty$. With this choice of $\de$, $K'=
\dE(e^{\nu_n^2(t_n-t_{n-1}) \sup\limits_{0\le u\le t_n} \al_u^1})
<\infty$ and thus:
$$\dE (|\zeta_M|) \le K\,K'\, \left( \dE \left( \langle x_M,x_M
    \rangle_{t_{n-1}}^{t_n}-(t_n^2-t_{n-1}^2)\right)^2 \right)^{1/2}\,\underset{M\to \infty}{\longrightarrow} 0$$
(by the law of large numbers for square-integrable independent
variables). Since $\arrowvert \frac{\cM_{t_n}}{\cM_{t_{n-1}}} \arrowvert \le 1$,
we also have $$\frac{\cM_{t_n}}{\cM_{t_{n-1}}} \,\zeta_M\,\overset{\dL^1}{\underset{M\to \infty}{\longrightarrow}} 0 \, .$$
Therefore
\begin{equation}
\label{eeeeee}
\dE(\frac{\cM_{t_n}}{\cM_{t_{n-1}}} \,\zeta_M \,|\,\cF_{t_{n-1}})
\,\overset{\dL^1}{\underset{M\to \infty}{\longrightarrow}} 0 \, .
\end{equation} 
In turn, by looking at (\ref{qwer}), this means that $$\dE(e^{i\nu_n(x_M(t_n)-x_M(t_{n-1}))}
\, |\, \cF_{t_{n-1}})\,\overset{\dL^1}{\underset{M\to
    \infty}{\longrightarrow}}\,
e^{\frac{\nu_n^2}{2}(t_n^2-t_{n-1}^2)}\, .$$ 
Now, plug this convergence and the induction hypothesis for $n-1$ into
(\ref{hypoa}) to get the result for $n$.

The same is true for $y_M$. To check that the finite-dimensionnal
distributions of $Z_M^{ij}$ have the right convergence, we would have
to prove that : 
$$\dE \bigg( \exp \Big( i \sum\limits_{i=1}^n \nu_i
    (x_M(t_i)-x_M(t_{i-1})) + \mu_i
    (y_M(t_i)-y_M(t_{i-1}))\Big)\bigg)$$ 
\begin{equation}
\label{wsx}
\underset{M \to  \infty}{\longrightarrow} 
\exp \left( \sum\limits_{i=1}^n
\frac{\nu_i^2+\mu_i^2}{2}\, (t_{i}^2-t_{i-1}^2)\right)\, .
\end{equation}
But since $ \langle x_M,y_M \rangle =0$, $$\cM_t=\exp\left(i(\nu_n x_M(t) + \mu_n y_M(t))-\frac{\nu_n^2}{2} \langle x_M,x_M \rangle_t -\frac{\mu_n^2}{2} \langle
  y_M,y_M \rangle_t\right)$$ is a martingale and the reasoning is
exactly the same as the previous one. 

Finally, let us look at the asymptotic independence. For the sake of
simplicity, let us take only two entries. Set for example
$x_M=Z^{11}_M$ and $y_M=\sqrt{2}\ree (Z^{12}_M)$. Then we have to
prove (\ref{wsx}) for our new $x_M,y_M$. Since $ \langle
x_M,y_M \rangle \not= 0$, $\cM_t$ previously defined is no more a
martingale. But $$\cN_t=\exp\left(i(\nu_n x_M(t) +  \mu_n
  y_M(t))-\frac{\nu_n^2}{2} \langle  x_M,x_M \rangle_t
  -\frac{\mu_n^2}{2} \langle  y_M,y_M \rangle_t - \nu_n \mu_n \langle
  x_M,y_M \rangle_t\right)$$ is a martingale and the fact that $\langle  x_M,y_M \rangle_t \overset{\dL^2}{\underset{M\to \infty}{\longrightarrow}} 0$ allows us to go along the same lines as before.

\end{proof}

\bibliography{bibiblio}
\bibliographystyle{amsalpha}

\begin{center}
  \hrule
\end{center}

\noindent
\textbf{E-mail}: \verb+Yan.Doumerc@math.ups-tlse.fr+.\\
\textbf{Postal Adress}: Laboratoire de Statistique et Probabilit\'es, U.M.R.
C.N.R.S. C5583, Universit\'e Paul Sabatier, 118 route de Narbonne, 31062
Toulouse CEDEX 4, France.

\end{document}